\documentclass[leqno]{article}
\usepackage{amsmath,amssymb,amsthm}
\usepackage[all]{xy}

\usepackage{amssymb} 
\usepackage{amsmath} 

\def\L{{\mathbf{L}}}
\def\N{{\mathbb{N}}}
\def\Z{{\mathbb{Z}}}
\def\R{{\mathbb{R}}}
\def\P{{\mathbb{P}}}
\def\H{{\text{H}}}

\DeclareMathOperator{\cohproj}{qgr}
\DeclareMathOperator{\proj}{proj}
\DeclareMathOperator{\coh}{coh}
\DeclareMathOperator{\Cone}{Cone}
\DeclareMathOperator{\Aut}{Aut}
\DeclareMathOperator{\Out}{Out}
\DeclareMathOperator{\QGr}{QGr}
\DeclareMathOperator{\Gr}{Gr}
\DeclareMathOperator{\gr}{gr}
\DeclareMathOperator{\qgr}{qgr}
\DeclareMathOperator{\Hom}{Hom}
\DeclareMathOperator{\Ext}{Ext}

\DeclareMathOperator{\DPic}{DPic}
\DeclareMathOperator{\End}{End}

\DeclareMathOperator{\Tor}{Tor}
\DeclareMathOperator{\tor}{tor}
\DeclareMathOperator{\Spaces}{Spaces}

\newcommand{\shift}{\stackrel{[1]}{\longrightarrow}}
\newcommand{\rmod}{\mathsf{mod}\text{-}A}

\newtheorem{theorem}{Theorem}[section]
\newtheorem{prop}[theorem]{Proposition}
\newtheorem{cor}[theorem]{Corollary}
\newtheorem{lemma}[theorem]{Lemma}
\newtheorem{define}[theorem]{Definition}
\newtheorem{rem}[theorem]{Remark}

\newtheorem{example}[theorem]{Example}

\title{A noncommutative version of Beilinson's Theorem}
\author{Hiroyuki Minamoto}

\begin{document}

\maketitle

\begin{abstract}

We prove that the category of representations of
the $N$-Kronecker quiver and
the category 
$\qgr\left(k\langle X_1,\dots,X_N \rangle/(\sum^N_{i=1}X_i^2)\right)$ 
which is considered as the category of coherent sheaves on
the noncommutative projective scheme associated to
$k\langle X_1,\dots,X_N \rangle/(\sum^N_{i=1}X_i^2)$
are derived equivalent.
This equivalence is easily proved by applying Orlov's Theorem,
on the other hand, 
in our proof, the quadratic relation $\sum_{i=1}^NX_i^2$
naturally arises from  Auslander-Reiten theory.

\end{abstract}
\setcounter{section}{-1}

\section{Introduction}$ $

In \cite{Be} 
Beilinson proved that 
the derived category of $\P^n$ is generated by 
the line bundles $\mathcal{O}_{\P^n},\dots,\mathcal{O}_{\P^n}(n)$.
This is equivalent to saying that 
the functor 
\[
\R\Hom(\bigoplus^{n}_{i=0}\mathcal{O}_{\P^n}(i),-):D^b(\coh \P^n)\longrightarrow D^b(\mathsf{mod}\text{-}B)
\]
is an equivalence of triangulated categories where 
$B=\End(\bigoplus^{n}_{i=0}\mathcal{O}_{\P^n}(i))$.
(This equivalence is also called Beilinson's Theorem.)

Although the projective space $\P^n_k$ is a basic object
in algebraic geometry,
the corresponding finite dimensional algebra $B$ is not elementary:

$B$ is the path algebra $k\overrightarrow{Q}_{n+1}$ 
of the Beilinson quiver $\overrightarrow{Q}_{n+1}$
\begin{equation*}
\begin{xymatrix}
{
0\ar@<-1ex>[r]_{x^{(0)}_{n+1}}\ar@{}[r]|-{:}\ar@<1ex>[r]^{x^{(0)}_1}&
1\ar@<-1ex>[r]_{x^{(1)}_{n+1}}\ar@{}[r]|-{:}\ar@<1ex>[r]^{x^{(1)}_1}&
2\ar@<-1ex>[r]_{x^{(2)}_{n+1}}\ar@{}[r]|-{:}\ar@<1ex>[r]^{x^{(2)}_1}&
\ar@{.}[r]&
\ar@<-1ex>[r]_{x^{(n-2)}_{n+1}}\ar@{}[r]|-{:}\ar@<1ex>[r]^{x^{(n-2)}_1}&
n-1
\ar@<-1ex>[r]_{\qquad x^{(n-1)}_{n+1}}
\ar@{}[r]|-{:}\ar@<1ex>[r]^{\qquad x^{(n-1)}_1}&n
}
\end{xymatrix}
\end{equation*}
with the relations
\[
x^{(l+1)}_ix^{(l)}_j=x^{(l+1)}_jx^{(l)}_i \quad 
\text{for}\,i,j\in\{1,\dots,n\}\,\text{and}\,l=0,\dots,n-2.
\]

In this paper we consider the path algebra $k\overrightarrow{\Omega}_N$ 
of $N$-Kronecker quiver $\overrightarrow{\Omega}_N$ 
\begin{equation*}
\begin{xymatrix}
{
0\,\ar@<-1ex>[r]_{x_N}\ar@{}[r]|-{:}\ar@<1ex>[r]^{x_1}&\,1
}
\end{xymatrix}
\end{equation*}
which is one of the simplest non-trivial finite dimensional algebra.
Our  main theorem states that
in noncommutative projective algebraic geometry,
there is an imaginary geometric object $\proj R$ which is 
derived equivalent to $\mathsf{mod}\text{-}k\overrightarrow{\Omega}_N$.

The main result of this paper is the following.
\begin{theorem}[Theorem \ref{main2}]\label{main}
Let $k$ be an algebraically closed field,
and $A=k\overrightarrow{\Omega}_N$, the path algebra of
the $N$-Kronecker quiver $\overrightarrow{\Omega}_N$.

Set $R=k\langle X_1,\dots,X_N \rangle/(\sum^N_{i=1}X_i^2)$.
Then for $N\geq 2$,
there exists a natural  equivalence of derived categories
\[
D^b(\rmod )\cong D^b(\qgr R),
\]
where $\rmod $ is the category of finite $A$-modules 
and 
$\qgr R$ is the category which is considered as 
the category of coherent sheaves 
on the noncommutative projective scheme $\proj R$ 
associated to the graded ring $R$.
\end{theorem}

This theorem is easily proved 
by applying Orlov's Theorem \cite[Theorem 2.5]{O}.
Orlov's Theorem gives 
a sufficient condition 
for a graded ring $R$ 
that 
$\proj R$ is derived equivalent to 
$\mathsf{mod}\text{-}k\overrightarrow{\Omega}_N$.
But in contrast, 
our theorem states that 
for the path algebra $k\overrightarrow{\Omega}_N$ 
there is a canonical choice of such $R$.

To prove Theorem \ref{main},
we use a method of noncommutative projective geometry 
and 
a result on the derived Picard group
by J. Miyachi and A. Yekutieli given in \cite{MY}.
They computed the derived Picard groups
of finite dimensional path algebras of  quivers
by using Happel's derived version of Auslander-Reiten theory,
from which
the quadratic relation $\sum_{i=1}^NX_i^2$
naturally arises.

In \cite{KR}, Kontsevich and Rosenberg constructed
the category  $\Spaces_k$ of ``noncommutative spaces,"
and showed that
in the category $\Spaces_k$,
the functor represented by
projective space $\P^n_k$ in (commutative) algebraic geometry
is represented by an object
$\text{N}\P^n_k$, the noncommutative projective space.
They also  proved that the category
$\coh \text{N}\P^n_k$ of coherent sheaves on  
$\text{N}\P^n_k$ is derived equivalent to
the category
$\mathsf{mod}\text{-}k\overrightarrow{\Omega}_{n+1}$.
Therefore it is natural to ask whether
the graded ring 
$k\langle X_1,\dots,X_{n+1}\rangle/(\sum_{i=1}^{n+1}X_i^2)$
is the coordinate ring of
$\text{N}\P^n_k$,i.e.,
whether  $\coh \text{N}\P^n_k$ is equivalent to
$\qgr k\langle X_1,\dots,X_{n+1}\rangle/(\sum_{i=1}^{n+1}X_i^2)$.

We  organize  the present  paper  as follows:
in Section 1
we introduce some definitions and results:
in Section 2 we prove Theorem \ref{main} (Theorem \ref{main2}).
\medskip

\textbf{Acknowledgment.}
I am grateful to I. Mori for his great help 
and comment on the first version of the paper.
I am also grateful to T. Abe
for his encouragement
 and
A. Moriwaki
for his hospitality and kindness.
I thank A. Takahashi for pointing out the paper of Orlov \cite{O}.
\section{Preliminaries}

In this section
we introduce some definitions and results.
Some of them are  about noncommutative projective schemes and
the others are about derived Picard groups.

\subsection{Noncommutative Projective Schemes}

\quad

This subsection is a summary of the paper \textup{\cite{Po}}
by A. Polishchuk.
We only treat $\N$-graded algebras
although $\Z$-algebras which is the more general notion
 are treated in \textup{\cite{Po}}.

Let $k$ be a field and
let $R=k\oplus R_1\oplus R_2\oplus\cdots$
be a connected graded coherent ring.
$\Gr R$ (resp. $\gr R$) denote the category of
graded right $R$-modules (resp. finitely presented graded right 
$R$-modules).
$\Tor R$ (resp. $\tor R$) denote the full subcategory of torsion modules
(resp. modules finite dimensional over $k$).
Note that $\Tor R$ and $\tor R$ are dense subcategories of
$\Gr R$ and $\gr R$ respectively, hence the quotient categories
$\QGr R=\Gr R/\Tor R$ and $\qgr R=\gr R/\tor R$ are abelian categories.

For graded a right $R$-module $M=\oplus_{n\in\Z} M_n$,
we denote  by $M(1)$ the $1$-degree shift of $M$.
i.e.,$M(1)_n=M_{n+1}$. 
The degree shift operator $(1): \gr R \longrightarrow \gr R$
induces the autoequivalence $(1)$ on $\qgr R$.
We denote by $\overline{R}$
the image in $\cohproj R$ of the regular module $R_R$.
The (coherent) \textit{noncommutative projective scheme} $\proj R$
associated to $R$ is the triple $(\qgr R, \overline{R},(1))$.
The autoequivalence $(1)$ is called the \textit{canonical polarization} on
$\proj R$.

In noncommutative projective geometry,
one thinks of $\qgr R$ 
as the category of coherent sheaves 
on a noncommutative projective scheme $\proj R$
associated to a graded ring $R$
($\cite{NPS,SvB,Po}$).

Let us consider 
a triple $(\mathcal{C},\mathcal{O},s)$ consisting of
a $k$-linear abelian category $\mathcal{C}$,
an object $\mathcal{O}$
and an autoequivalence $s$ on $\mathcal{C}$.
For $\mathcal{F}\in \mathcal{C}$, we define
\[
\Gamma_{*}(\mathcal{F})=
\bigoplus_{n\geq 0}\Hom_{\mathcal{C}}(\mathcal{O},\mathcal{F}(n) ),
\]
where $\mathcal{F}(n)=s^n\mathcal{F}$,
and we set
\[
R=\Gamma_{*}(\mathcal{C},\mathcal{O},s)=\Gamma_{*}(\mathcal{O}).
\]

Multiplication is defined as follows:
If $x\in \Hom_{\mathcal{C}}(\mathcal{O},\mathcal{F}(l))$,
$b \in \Hom_{\mathcal{C}}(\mathcal{O},\mathcal{O}(m))$ and
$a \in \Hom_{\mathcal{C}}(\mathcal{O},\mathcal{O}(n))$ then
\[
x\cdot a=s^n(x)\circ a
\quad
\text{and}
\quad
a\cdot b = s^m(a)\circ b.
\]
With this law of composition,
$\Gamma_{*}(\mathcal{F})$ becomes a graded right module
 over  the graded algebra $R$ over $k$.

\begin{define}
[\textup{\cite[Section 4.2]{NPS},\cite[Section 2]{Po}}]
\label{ample}
Let $(\mathcal{C},\mathcal{O},s)$ be a triple as above.
Then the pair $(\mathcal{O},s)$ is called ample 
if the following conditions hold:

\begin{itemize}
\item[(1)]
For every object $\mathcal{F}\in\mathcal{C}$,
there are positive integers $l_1,\dots,l_p$ and
an epimorphism $\oplus_{i=1}^p\mathcal{O}(-l_i)\longrightarrow\mathcal{F}$.

\item[(2)]
For every epimorphism $f: \mathcal{F}\longrightarrow\mathcal{G}$,
there exists an integer $n_0$ such that for every $n\geq n_0$
the induced map
$\Hom_{\mathcal{C}}(\mathcal{O},\mathcal{F}(n))\longrightarrow
\Hom_{\mathcal{C}}(\mathcal{O},\mathcal{G}(n))$
is surjective.

\end{itemize}
\end{define}

Let $\pi: \Gr R\longrightarrow \QGr R$ be the quotient functor.
Set $\overline{\Gamma}_*=\pi\circ\Gamma_*$.

\begin{theorem}
[\textup{\cite[Theorem 2.4]{Po}}]\label{Serre}
Let $(\mathcal{C},\mathcal{O},s)$ be a triple as above.
If $(\mathcal{O},s)$ is ample,
then the graded ring
$R=\Gamma_{*}(\mathcal{C},\mathcal{O},s)$ is coherent,
$\Gamma_*(\mathcal{F})$ is finitely presented $R$-module
for $\mathcal{F}\in\mathcal{C}$
and the functor
$\overline{\Gamma}_*:\mathcal{C}\longrightarrow \cohproj R$
induces an equivalence of triples 
between 
$(\mathcal{C},\mathcal{O},s)$ and 
$\proj R =(\qgr R,\overline{R},(1))$,
i.e.,
\begin{itemize}
\item[] $\overline{\Gamma}_*:\mathcal{C}\stackrel{\sim}{\longrightarrow}\cohproj R$
is an equivalence of categories,
\item[] $\overline{\Gamma}_*(\mathcal{O})\cong \overline{R}$,
and $\overline{\Gamma}_*\circ s = (1)\circ\overline{\Gamma}_*$.
\end{itemize}
\end{theorem}

By Beilinson's Theorem 
$\coh \P^1$ is derived equivalent to
$\mathsf{mod}\text{-}k\overrightarrow{\Omega}_2$.
Therefore by  Theorem \ref{main},
$\coh\P^1$ and
$\qgr \left(\frac{k\langle X_1,X_2 \rangle}{(X_1^2+X_2^2)}\right)$ 
are derived equivalent.
However, it turns out that these two categories are equivalent.

\begin{example}
[\textup{\cite[Section 3]{SvB}}] 
Let $\sigma : \P^1\longrightarrow \P^1,\,[a:b]\mapsto [b:-a]$
be the automorphism of $\P^1$ and
set $s=\sigma_*(-\otimes_{\mathcal{O}_{\P^1}}\mathcal{O}_{\P^1}(1))$.
Then $s$ is ample on the triple
$(\coh (\P^1),\mathcal{O}_{\P^1},s)$ and
\[
\Gamma_{*}(\P^1,\mathcal{O}_{\P^1},s)\cong
\frac{k\langle X_1,X_2\rangle}{(X_1^2+X_2^2)}\]

Traditionally the graded ring $\Gamma_{*}(\P^1,\mathcal{O}_{\P^1},s)$
is denoted by $B(\P^1,\mathcal{O}_{\P^1}(1),\sigma)$ and
is called the twisted homogeneous coordinate ring
associated to the pair $(\mathcal{O}_{\P^1}(1),\sigma)$.

Hence by  Theorem \ref{Serre}, 
the triples 
$(\coh\P^1,\mathcal{O}_{\P^1},s)$ and 
$\proj \left(\frac{k\langle X_1,X_2\rangle}{(X_1^2+X_2^2)}\right)$ 
are equivalent.
In particular there is an equivalence
\[
\overline{\Gamma}_* : 
\coh \P^1 \stackrel{\sim}{\longrightarrow}
\qgr \left(\frac{k\langle X_1,X_2\rangle}{(X_1^2+X_2^2)}\right) .
\]
\end{example}

Moreover, 
from the proof of Theorem \ref{main} in Section 2,
it is easy to see that 
there is a commutative diagram 
\begin{equation*}
\begin{xymatrix}
{
D^b(\coh\P^1)\ar[r]^{\sim}&
D^b(\mathsf{mod}\text{-}k\overrightarrow{\Omega}_2)&
D^b\left(\qgr\left(\frac{k\langle X_1,X_2\rangle}{(X_1^2+X_2^2)}\right)\right)
{\ar[l]_{\sim\quad}}\\
\coh\P^1\ar@{^{(}->}[u] \ar[rr]^{\sim}_{\overline{\Gamma}_*} & &
\qgr\left(\frac{k\langle X_1,X_2\rangle}{(X_1^2+X_2^2)}\right)
\ar@{^{(}->}[u]
}
\end{xymatrix}
\end{equation*}
where the vertical arrows are natural inclusions.

\begin{rem}
The triples 
$(\coh\P^1,\mathcal{O}_{\P^1},-\otimes_{\P^1}\mathcal{O}_{\P^1}(1))$ 
and 
$\proj(k[X_1,X_2])$ are 
equivalent 
by  Serre's Theorem $\textup{\cite[Proposition II.5.15]{AG}}$.
Hence 
the categories 
$\qgr (k[X_1,X_2])$ and  
$\qgr\left(\frac{k\langle X_1,X_2\rangle}{(X_1^2+X_2^2)}\right)$ 
are equivalent.

Moreover these two (noncommutative) projective schemes are equivalent 
at the level of affine covers.
There is an  equivalence 
\[
\gr (k[X_1,X_2])
\stackrel{\sim}{\longrightarrow} 
\gr\left(\frac{k\langle X_1,X_2\rangle}{(X_1^2+X_2^2)}\right)
\]
which sends $k[X_1,X_2]$ to 
$\frac{k\langle X_1,X_2\rangle}{(X_1^2+X_2^2)}$. 
As $\tor$ is 
the full subcategory of $\gr$ 
consisting of objects of finite length,
from the above equivalence of affine covers 
we can also obtain the equivalence between 
$\qgr(k[X_1,X_2])$ and 
$\qgr\left(\frac{k\langle X_1,X_2\rangle}{(X_1^2+X_2^2)}\right) $.
\textup{(See \cite[Section 3]{Z})}. 
\end{rem}

\subsection{Derived Picard Groups of Finite Dimensional Hereditary Algebras}

$ $

This subsection is a summary of \textup{\cite{MY}}
by J. Miyachi and A. Yekutieli.

First note that we  work with right modules,
hence the definition of Riedtmann quiver differs from that of
\cite{Ha} and \cite{MY}. (They work with left modules.)

Let $A$ be a finite dimensional $k$-algebra.
A complex $T\in D^b(\mathsf{Mod}\text{-}A\otimes_kA^{\text{op}})$
is called a \textit{two-sided tilting complex}
if there exists another complex
$T^{\vee}\in D^b(\mathsf{Mod}\text{-}A\otimes_kA^{\text{op}})$
such that
$T\otimes^\L_{A}T^{\vee}\cong T^{\vee}\otimes^{\L}_AT\cong A$.
The \textit{derived Picard group} of $A$ (relative to $k$) is
\[
\DPic_k(A):= \{\text{two-sided tilting complexes}\}/\text{isomorphisms}
\]
with the identity element $A$,
product $(T_1,T_2)\mapsto T_1\otimes^\L_AT_2$ and
inverse $T\mapsto T^\vee := \R\Hom_A(T,A)$.

A tilting complex induces an equivalence of triangulated categories,
\[
-\otimes^\L_AT: D^b(\mathsf{Mod}\text{-}A)\longrightarrow
D^b(\mathsf{Mod}\text{-}A).
\]
See \cite{Y},\cite{MY} for more details.

Let
$\overrightarrow{\triangle}=
(\overrightarrow{\triangle}_0,\overrightarrow{\triangle}_1)$
be a finite quiver.
The Riedtmann quiver $\overrightarrow{\Z}\overrightarrow{\triangle}$
of $\overrightarrow{\triangle}$ is defined as follows :
The set of vertices $(\overrightarrow{\Z}\overrightarrow{\triangle})_0$
is given by $\Z\times \overrightarrow{\triangle}_0$.
Given an arrow $a:x\rightarrow y$ in $\overrightarrow{\triangle}$ 
and an integer $n$,
there are  arrows $(n,a):(n,y)\rightarrow (n,x)$ and
$(n,a)^*:(n,x)\rightarrow (n+1,y)$.
The translation $\tau$ and
the polarization $\mu$ on
$\overrightarrow{\Z}\overrightarrow{\triangle}$ are defined by
$\tau(n,x)=(n-1,x)$
and $\mu(n,a)=(n-1,a)^*,\,\mu((n,a)^*)=(n,a)$.
For an arrow $\alpha: \xi\rightarrow \eta $ in
$\overrightarrow{\Z}\overrightarrow{\triangle}$ ,
we have $\mu(\alpha):\tau(\eta)\rightarrow \xi$.

The \textit{path category}
$k\langle\overrightarrow{\Z}\overrightarrow{\triangle}\rangle$
(in the sense of \cite{MY}) is 
the  category
whose set of objects is $(\overrightarrow{\Z}\overrightarrow{\triangle})_0$,
morphisms are generated by identities and the arrows and the only relations
arise from incomposability of paths.

Let $\eta$ be a vertex of  $\overrightarrow{\Z}\overrightarrow{\triangle}$,
let $\xi_1,\cdots,\xi_p$ be complete representatives of
the set
$\{\xi\mid \text{there is an arrow}\,\,\alpha:\xi\rightarrow \eta\}$
and let $\{\alpha_{ij}\}_{j=1}^{d_i}$ be the set of arrows
from $\xi_i$ to $\eta$.
The \textit{mesh} ending at $\eta$ is
the subquiver of $\overrightarrow{\Z}\overrightarrow{\triangle}$
with vertex
$\{\eta,\mu(\eta),\xi_1,\dots,\xi_p\}$
and arrows
$\{\alpha_{ij},\mu(\alpha_{ij})\mid i=1,\dots,p,\,j=1,\dots,d_i\}$.
The \textit{mesh ideal} in the path category
$k\langle\overrightarrow{\Z}\overrightarrow{\triangle}\rangle$
is the ideal
generated by the elements
\[
\sum_{i=1}^p\sum_{j=1}^{d_i}\alpha_{ij}\circ\mu(\alpha_{ij})\,
\in
\Hom_{ k\langle\overrightarrow{\Z}\overrightarrow{\triangle}\rangle}
(\mu(\eta),\eta).
\]
The \textit{mesh category}
$k\langle\overrightarrow{\Z}\overrightarrow{\triangle},I_m\rangle$
is defined as the quotient category of the path category modulo mesh ideal.

Let $A=k\overrightarrow{\triangle}$ be the path algebra of
$\overrightarrow{\triangle}$ and
let $\rmod$ denote the category of finite right $A$-modules.
The $k$-dual $A^*=\Hom_k(A,k)$ of $A$ is a two-sided tilting complex.
We write $\tau_A\in\DPic_k(A)$ for the element represented by $A^*[-1]$.
Let us agree that  $\tau_A$ also  denotes the autequivalence
$-\otimes_A^{\L}A^*[-1]$ of $D^b(\rmod)$.

Let $P_x$ be the indecomposable projective right $A$-module
corresponding to a vertex
$x\in\overrightarrow{\triangle}_0$.
Define $\mathsf{B}\subset D^b(\rmod)$ to be the full subcategory
with objects
$\{\tau_{A}^n P_x\mid x\in\overrightarrow{\triangle}_0,n\in\Z\}$.

The following theorem which is important in \cite{MY} 
is also important for us.

\begin{theorem}
[\textup{\cite[Theorem 2.6]{MY}}]\label{mesh}
There is a $k$-linear equivalence
\[
G:  k\langle\overrightarrow{\Z}\overrightarrow{\triangle},I_m\rangle
\longrightarrow \mathsf{B},\quad
G(n,x)=\tau_A^{-n}P_x.
\]
\end{theorem}

The equivalence $G$ sends
the mesh ending at $\eta=(n,y)$ to the exact triangle
called
\textit{Auslander-Reiten triangle}
\begin{equation}\label{A-R-tr}
\tau_A(G(\eta))=
G(\tau(\eta))\stackrel{\oplus G(\beta_{ij})}{\longrightarrow}
\bigoplus_{i=1}^{p}\bigoplus_{j=1}^{d_i}G(\xi_i)
\stackrel{\oplus G(\alpha_{ij})}{\longrightarrow}
G(\eta)\shift
\end{equation}
where $\beta_{ij}=\mu(\alpha_{ij})$.
By the definition of the path category
the sets of arrows $\{\alpha_{ij}\}^{d_i}_{j=1}$ and
$\{\beta_{ij}\}^{d_i}_{j=1}$ are basis of 
$\Hom_{k\langle\overrightarrow{\Z}\overrightarrow{\triangle},I_m\rangle}
(\xi_i,\eta)$
and
$\Hom_{k\langle\overrightarrow{\Z}\overrightarrow{\triangle},I_m\rangle}
(\tau(\eta),\xi_i)$ respectively.
Fix isomorphisms
\begin{equation*}
\begin{split}
\Hom_{k\langle\overrightarrow{\Z}\overrightarrow{\triangle},I_m\rangle}
(\xi_i,\eta)&\cong
\Hom_{k\langle\overrightarrow{\Z}\overrightarrow{\triangle},I_m\rangle}
(\tau(\eta),\xi_i)^*\\
 \alpha_{ij}&\mapsto \beta^*_{ij}
\end{split}
\end{equation*}
where $\{\beta_{ij}^*\}_{j=1}^{d_i}$ is the dual basis of
$\{\beta_{ij}\}_{i=1}^{d_i}$.
Set
$V_i=
\Hom_{k\langle\overrightarrow{\Z}\overrightarrow{\triangle},I_m\rangle}
(\xi_i,\eta)\cong
\Hom_{k\langle\overrightarrow{\Z}\overrightarrow{\triangle},I_m\rangle}
(\tau(\eta),\xi_i)^*$.
Then there are canonical morphisms
\begin{equation}\label{cm-A-R}
\varphi_i : V_i\otimes_k G(\xi_i)\longrightarrow G(\eta),\quad
\psi_i : \tau_A(G(\eta)) \longrightarrow V_i\otimes_k G(\xi_i).
\end{equation}
Then 
the Auslander-Reiten triangle (\ref{A-R-tr}) has the form 
\begin{equation}\label{A-R-tr2}
\tau_A(G(\eta))\stackrel{\oplus \psi_i}{\longrightarrow}
\bigoplus_{i=1}^{p}V_i\otimes_kG(\xi_i)
\stackrel{\oplus \varphi_i}{\longrightarrow} G(\eta) \shift .
\end{equation}

Let $\Aut((\overrightarrow{\Z}\overrightarrow{\triangle})_0)$ be
the permutation group of the vertex set
$(\overrightarrow{\Z}\overrightarrow{\triangle})_0$
and let
$\Aut((\overrightarrow{\Z}\overrightarrow{\triangle})_0;d)
^{\langle\tau\rangle}$ be the subgroup of permutations which preserve
arrow-multiplicities and 
which commute with $\tau$, 
namely
\begin{align*}
&\Aut((\overrightarrow{\Z}\overrightarrow{\triangle})_0;d)
^{\langle\tau\rangle}=   \\
&\{\pi\in\Aut(\overrightarrow{\Z}\overrightarrow{\triangle})_0\mid
d(x,y)=d(\pi(x),\pi(y))
\,\text{for all}\,
x,y\in (\overrightarrow{\Z}\overrightarrow{\triangle})_0 ,
\,\text{and}\,
\pi\tau=\tau\pi\}
\end{align*}
where $d(x,y)$ denotes the arrow-multiplicity from $x$ to $y$.

Now we can state one of the main result of \cite{MY},
which shall be used to prove 
our main theorem.

\begin{theorem}
[\textup{\cite[Theorem 3.8]{MY}}]\label{DPic-th}
Let $\overrightarrow{\triangle}$ be a finite quiver without oriented cycles
and let $A=k\overrightarrow{\triangle}$ be the path algebra of
$\overrightarrow{\triangle}$ over an algebraically closed field $k$.
If $A$ has infinite representation type
then there is an isomorphism of groups
\[
\DPic_k{A}\cong
(\Aut((\overrightarrow{\Z}\overrightarrow{\triangle})_0;d)\ltimes
\Out_k^0(A))\times \Z
\]
where $\Out_k^{0}(A)$ denotes the identity component of
the group of outer automorphisms $\Out_k(A)$.
\end{theorem}

\begin{rem}
In \textup{\cite{MY}},
the case of the finite representation type is also computed.
\end{rem}

\section{Proof of Theorem 0.1}

$ $

From now on we assume that
$k$ is an algebraically closed field and $N\geq 2$.

Let $A$ be the path algebra of
the  $N$-Kronecker quiver $\overrightarrow{\Omega}_N$ shown in introduction.
We denote by $P_0,P_1$ the indecomposable projective right $A$-modules
associated to vertices $0,1$ of $\overrightarrow{\Omega}_N$,
and $S_0,S_1$ the simple $A$ modules
associated to vertices $0,1$ of $\overrightarrow{\Omega}_N$.
Then $S_1=P_1$ and $P_0$ is the projective cover of $S_0$.

The Riedetmann quiver $\overrightarrow{\Z}\overrightarrow{\Omega}_N$
of the $N$-Kronecker quiver $\overrightarrow{\Omega}_N$ is
shown in Figure 1 below.
(In Figure 1 vertices $(-n,1)$ and $(-n,0)$ are
replaced by $\tau^n_AP_1$ and $\tau^n_AP_0$ respectively.)

\begin{figure}[htbp]\label{R-q}
\begin{equation*}
\begin{xymatrix}
{
\tau^{n+1}P_1\ar@<-1ex>[dr]_{x_N}\ar@{}[dr]|-{:}\ar@<1ex>[dr]^{x_1}&&
\tau^{n}P_1\ar@<-1ex>[dr]_{x_N}\ar@{}[dr]|-{:}\ar@<1ex>[dr]^{x_1}&&
\tau^{(n-1)}P_1\\
&\tau^{n+1}P_0\ar@<-1ex>[ur]_{x_N}\ar@{}[ur]|-{:}\ar@<1ex>[ur]^{x_1}&&
\tau^{n}P_0\ar@<-1ex>[ur]_{x_N}\ar@{}[ur]|-{:}\ar@<1ex>[ur]^{x_1}
}
\end{xymatrix}
\end{equation*}
\caption{Riedetmann quiver
$\protect\overrightarrow{\Z}\protect\overrightarrow{\Omega}_N$ }
\end{figure}

It is easy to see that
\[
\Aut((\overrightarrow{\Z}\overrightarrow{\Omega}_N)_0;d)\cong \Z
\]
and there is a generator $\rho$ of 
$\Aut((\overrightarrow{\Z}\overrightarrow{\Omega}_N)_0;d)$
 such that $\rho(0,1)=(0,0),\rho(0,0)=(1,1)$.
This $\rho$ satisfies the relation $\rho^{-2}=\tau$.
By Theorem \ref{DPic-th}
there exists a two-sided tilting complex $\rho_A$ such that
$\rho_A^{-2}\cong \tau_A$ and $\rho_A P_1\cong P_0$.

From now on we write $\rho=\rho_A$ and $\tau=\tau_A$.

For $M^\cdot\in D^b(\rmod)$  we use the following notation
\[
\rho^{n}M^\cdot=
M^\cdot\otimes^\L_A
\overbrace{\rho\otimes^\L_A\cdots\otimes^\L_A\rho}^n.
\]

\begin{figure}[htbp]
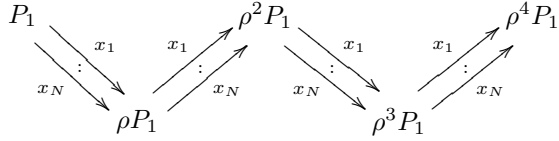
\label{R-q2}
\begin{equation*}
\begin{xymatrix}
{
P_1\ar@<-1ex>[dr]_{x_N}\ar@{}[dr]|-{:}\ar@<1ex>[dr]^{x_1}&&
\rho^2P_1\ar@<-1ex>[dr]_{x_N}\ar@{}[dr]|-{:}\ar@<1ex>[dr]^{x_1}&&
\rho^4P_1\\
&\rho P_1\ar@<-1ex>[ur]_{x_N}\ar@{}[ur]|-{:}\ar@<1ex>[ur]^{x_1}&&
\rho^3 P_1\ar@<-1ex>[ur]_{x_N}\ar@{}[ur]|-{:}\ar@<1ex>[ur]^{x_1}
}
\end{xymatrix}
\end{equation*}
\caption{Riedetmann quiver
$\protect\overrightarrow{\Z}\protect\overrightarrow{\Omega}_N$}
\end{figure}

Since $P_0\cong\rho P_1$ and $\rho^2\cong \tau^{-1}$ it follows that
$\tau^{-n}P_1\cong\rho^{2n}P_1,\,\tau^{-n}P_0\cong\rho^{2n+1}P_1$,
the Riedetmann quiver
$\overrightarrow{\Z}\overrightarrow{\Omega}_N$ has the 
form shown in Figure \ref{R-q2}.

Set $V=\Hom(P_1,\rho P_1)\cong\Hom(P_1,P_0)$.
The set of arrows $x_1,\dots,x_N$ is a basis of $V$.
The graded vector space 
$\bigoplus_{n\geq 0}\Hom(P_1,\rho^n P_1)$
has a graded algebra structure as in the same way of Section 1.1.

By  the mesh relation and Theorem \ref{mesh}
we have the following proposition.

\begin{prop}\label{coordinate}
There is an isomrphism of graded algebras
\[
\gamma:\bigoplus_{n\geq 0}\Hom(P_1,\rho^n P_1) \longrightarrow k\langle X_1,\dots,X_N\rangle/(\sum^N_{i=1}X_i^2)
\]
 which sends $x_i\in V$ to $X_i$ for $i=1,\dots,N$.
\end{prop}

Fix isomorphisms
$V\cong
\Hom_{D^b(\rmod)}(\rho^{n}P_1,\rho^{n+1}P_1)$ for $n\geq 0$ and
identify
$V$ with its dual $V^*$ by sending
$x_i\mapsto x_i^*$.
Then the canonical morphisms $\varphi_i$ and $\psi_i$
in (\ref{cm-A-R}) of Section 1.2
have  the following forms:
\begin{equation*}
\begin{split}
&\varphi: V\otimes_k\rho^nP_1\longrightarrow \rho^{n+1}P_1,\,
v\otimes p_n \mapsto v(p_n),\\
&\psi: \rho^{n-1}P_1\longrightarrow V\otimes \rho^nP_1,\,
p_{n-1}\mapsto \sum_{i=1}^{N}x_i\otimes x_i(p_{n-1}).
\end{split}
\end{equation*}

The Auslander-Reiten triangle (\ref{A-R-tr2}) has
the form 
\begin{equation*}
(\mathcal{S}_n)\qquad
 \rho^{n-1}P_1
\stackrel{\psi}{\longrightarrow}
V\otimes_k\rho^{n}P_1
\stackrel{\varphi}{\longrightarrow}\rho^{n+1}P_1\shift.
\end{equation*}

Set $R=k\langle X_1,\dots,X_N\rangle/(\sum^N_{i=1}X_i^2)$.
We denote by $R_n$ the degree $n$ component of the graded ring
$R=\bigoplus_{n\geq 0}R_{n}=k\langle X_1,\dots,X_N\rangle/(\sum^N_{i=1}X_i^2)$.
The following lemma is easy to show.

\begin{lemma}\label{ex-seq}
For each $n \geq 1$,
there is a following short exact sequence of $k$-vector spaces 
\begin{equation}\label{ex-seq-R}
0\longrightarrow R_{n-1}
\stackrel{\Psi}{\longrightarrow }R_1\otimes_k R_{n}
\stackrel{\Phi}{\longrightarrow} R_{n+1}
\longrightarrow 0,
\end{equation}
where the second and third morphisms are defined as follows:
\begin{equation*}
\begin{split}
\Phi : &R_1\otimes_k R_{n}\ni f\otimes_k r \mapsto fr \in R_{n+1}\\
\Psi : &R_{n-1}\ni r_{n-1}\mapsto
\sum_{i=1}^NX_i\otimes_k X_ir_{n-1}\in R_1\otimes_k R_n.
\end{split}
\end{equation*}
.
\end{lemma}

Let 
$\gamma_n:\Hom(P_1,\rho^nP_1)\stackrel{\sim}{\longrightarrow} R_n$ 
be the degree $n$ part of 
the isomorphism $\gamma$ of 
Proposition \ref{coordinate}.
Then we have the commutative diagram 
\begin{equation*}
\begin{xymatrix}
{
&\Hom(P_1,\rho^{n-1}P_1)\ar[rr]^{\Hom(P_1,\psi)}\ar[d]_{\gamma_{n-1}}^{\wr}& &
V\otimes_k\Hom(P_1,\rho^{n}P_1)\ar[d]_{\gamma_{1}\otimes\gamma_{n}}^{\wr}\\
0\ar[r]& R_{n-1}\ar[rr]_{\Psi}& &R_1\otimes_k R_{n} 
}
\end{xymatrix}
\end{equation*}
where 
the top row is obtained by applying $\Hom(P_1,-)$ 
to the first two terms of
the exact triangle $(\mathcal{S}_n)$ 
and the bottom row is 
the first three terms of the exact sequence (\ref{ex-seq-R}) 
of Lemma \ref{ex-seq}.
Note, 
since two vertical arrows in the diagram are isomorphisms, 
the injectivity of $\Psi$  
implies that $\Hom(P_1,\psi)$ is injective.

\begin{prop}\label{pure}
$\rho^n=\rho \otimes^{\L}_A\cdots\otimes^{\L}_A\rho$
is a pure module for $n\geq 0$,
\textup{(}
i.e.,$\H^{i}(\rho^n)=0$ for all $i\neq 0$
\textup{)}
.
\end{prop}

\begin{proof}
Since $\rho^n \cong \rho^n P_1\oplus \rho^{n+1} P_1$ and 
$\rho P_1\cong P_0$ are pure modules,
it suffices to show that $\rho^{n+1} P_1$ is pure for $n\geq 1$.
By induction on $n$ 
we may assume that $\rho^{n-1}P_1$ and $\rho^nP_1$ are pure.
Then
the long cohomology sequence of 
the Auslander-Reiten triangle 
$(\mathcal{S}_n)$ implies
that $\H^{i}(\rho^{n+1}P_1)= 0$ for
$i\neq -1,0$ .

Note that 
$\H^{-1}(\rho^{n+1})=0$ if and only if 
\[
\Ext^{-1}(A,\rho^{n+1}P_1)\cong\Hom(A,\H^{-1}(\rho^{n+1}P_1))=0.
\]
Since 
$A \cong P_1 \oplus \rho P_1$ and 
$\Ext^{-1}(P_1,\rho^{n}P_1)\cong \Hom(P_1,\H^{-1}(\rho^{n}P_1))=0$ 
by the induction hypothesis,
we have 
\[
\Ext^{-1}(A,\rho^{n+1}P_1)\cong \Ext^{-1}(P_1,\rho^{n+1}P_1)\oplus\Ext^{-1}(P_1,\rho^{n}P_1)
\cong \Ext^{-1}(P_1,\rho^{n+1}P_1).
\]
Thus we have only to show 
$\Ext^{-1}(P_1,\rho^{n+1}P_1)\cong \Hom(P_1,\H^{-1}(\rho^{n+1}P_1))= 0 $.

There is an exact sequence of $A-$modules
\begin{equation*}
0\rightarrow
\H^{-1}(\rho^{n+1}P_1)\longrightarrow
 \rho^{n-1}P_1
\stackrel{\psi}{\longrightarrow}
V\otimes_k\rho^{n}P_1
\longrightarrow\H^{0}(\rho^{n+1}P_1)\rightarrow 0.
\end{equation*}
Applying the exact functor $\Hom(P_1,-)$ 
to the first four terms of this sequence 
yields the following exact sequence
\begin{multline*}
0\longrightarrow \Hom(P_1,\H^{-1}(\rho^{n+1}P_1))\longrightarrow \\
\Hom(P_1,\rho^{n-1}P_1)
\stackrel{\Hom(P_1,\psi)}{\longrightarrow}
V\otimes_k\Hom(P_1,\rho^{n}P_1).
\end{multline*}
The injectivity of $\Hom(P_1,\psi)$ 
implies $\Hom(P_1,\H^{-1}(\rho^{n+1}P_1))=0$.
\end{proof}

\begin{cor}\label{pure-cor}
$ $

\begin{itemize}
\item[(1)]
$\tau^{n+1}[1]$ is a pure module for every $n\geq 0$.

\item[(2)]
$\rho^{-n} S_0\cong \rho^{-(n+1)}P_1[1]$ is 
a pure module for every $n\geq 0$.
\end{itemize}
\end{cor}

\begin{proof}
$(1)$
We have 
$\tau^{n+1}[1]\cong\tau^{n}\otimes^{\L}_A A^*
\cong\R\Hom(\tau^{-n},A^*)\cong (\tau^{-n})^*$, 
using  the property of 
the Serre functor $-\otimes_A^\L{A^*}$ 
for the third isomorphism.
By Proposition \ref{pure}
$\tau^{-n}\cong \rho^{2n}$ is pure for $n\geq 0$.
Hence $\tau^{n+1}[1]$ is also pure for $n\geq 1$

$(2)$
By the direct calculation we have
\[
S_0\cong P_0\otimes_A^\L A^*\cong (\tau P_0)[1] \cong (\rho^{-1}P_1)[1].
\]
Therefore 
$\rho^{-n}S_0\cong \rho^{-(n+1)}P_1[1]$. 
$\rho^{-2n}S_0\cong (\tau^{n+1}P_0)[1] \cong P_0\otimes_A \tau^{n+1}[1]$ and
$\rho^{-(2n+1)}S_0\cong (\tau^{n+1}P_1)[1]\cong P_1\otimes_A\tau^{n+1}[1]$ are pure by (1).
\end{proof}

\begin{rem}
In this paper we assume $N\geq 2$.
If $N=1$ then 
it is known that $\tau^3\cong A[-2]$ (\textup{\cite[Theorem 4.1]{MY}}).
Hence the above corollary fails.
\end{rem}

The following lemma 
shall be used in the proof of the Lemma \ref{t-st-lem}.

\begin{lemma}\label{surj-lem}
The canonical morphism 
\[
c:\Hom(\rho^{-(n+1)}S_0,\rho^{-n}S_0)\otimes_k\rho^{-(n+1)}S_0
\longrightarrow\rho^{-n}S_0
\]
is surjective.
\end{lemma}

\begin{proof}
The canonical morphism 
$c$
 is the 1-shift of 
the last two terms of
the Auslander-Reiten triangle 
\begin{equation*}
(\mathcal{S}_{-(n+2)})\qquad
 \rho^{-(n+3)}P_1
\stackrel{\psi}{\longrightarrow}
V\otimes_k\rho^{-(n+2)}P_1
\stackrel{\varphi}{\longrightarrow}\rho^{-(n+1)}P_1\shift.
\end{equation*}
under the isomorphisms $\rho^{-n} S_0\cong \rho^{-(n+1)}P_1[1]$ 
and $V\cong\Hom(\rho^{-(n+1)}S_0,\rho^{-n}S_0)$.
Therefore the cokernel of $c$ 
is $\H^{0}(\rho^{-(n+3)}P_1[2])=0$ 
by Corollary \ref{pure-cor}.(2). 
\end{proof}

Let us review the definition of a t-structure.
A t-structure in a triangulated category $\mathcal{T}$ 
is a pair $(\mathcal{T}^{\leq 0},\mathcal{T}^{\geq 0})$ 
of full subcategories which are closed under isomorphisms 
such that, 
for 
$\mathcal{T}^{\geq n}:= (\mathcal{T}^{\geq 0})[-n],n\in\Z$
and  
$\mathcal{T}^{\leq n}:= (\mathcal{T}^{\leq 0})[-n],n\in\Z$ 
the following three conditions are satisfied:
\begin{itemize}
\item[(t1)]
 For $X\in \mathcal{T}^{\leq 0}$ and $Y\in\mathcal{T}^{\geq 0}$ 
we have that $\Hom_{\mathcal{T}}(X,Y)=0$. 

\item[(t2)]
$\mathcal{T}^{\leq 0}\subset\mathcal{T}^{\leq 1}$ and 
$\mathcal{T}^{\geq 1}\subset\mathcal{T}^{\geq 0}$.

\item[(t3)]
For $X\in \mathcal{T}$ 
there is an exact triangle
$A\rightarrow X\rightarrow B\shift$
such that 
$A\in \mathcal{T}^{\leq 0}$ and 
$B\in \mathcal{T}^{\geq 1}$. 
\end{itemize}

The full subcategory 
$\mathcal{H}:=\mathcal{T}^{\geq 0}\cap\mathcal{T}^{\leq 0}$ 
is called the \textit{heart} of the t-structure 
$(\mathcal{T}^{\geq 0},\mathcal{T}^{\leq 0})$.
The important fact is that 
the heart $\mathcal{H}$ of a t-structure 
is an abelian category. (See \cite[IV.4]{GM})

Let $X$ be a projective variety over $k$ 
and $\mathcal{O}_X(1)$ be an ample line bundle on $X$.
Let $D^{\geq 0}(\coh X)$ (resp. $D^{\leq 0}(\coh X)$) be 
the full subcategory of $D^{b}(\coh X)$ 
with objects $\mathcal{F}^\cdot$ such that 
$\H^{i}(\mathcal{F}^\cdot)=0$ for $i < 0$ (resp. $i>0$).
Then the pair $(D^{\geq 0}(\coh X),D^{\leq 0}(\coh X))$ is 
a t-structure in $D^b(\coh X)$ and 
its heart is $\coh X$.
By Serre's vanishing theorem \textup{\cite[Propsition III.5.3]{AG}},
$D^{\geq 0}(\coh X)$ (resp. $D^{\leq 0}(\coh X)$) 
is a full subcategory 
with objects $\mathcal{F}^\cdot$ which satisfy 
\begin{equation*}
\begin{split}
&\R\Hom^{\cdot}(\mathcal{O}_X,\mathcal{F}^\cdot(n))
\in D^{\geq 0}(k\text{-}vect)
\quad
\text{for $n \gg 0$}\\
(\text{resp.}\quad
&\R\Hom^{\cdot}(\mathcal{O}_X,\mathcal{F}^\cdot(n))
\in D^{\leq 0}(k\text{-}vect)
\quad
\text{for $n \gg 0$})
\end{split}
\end{equation*}
where we denote 
\[
\mathcal{F}^\cdot(n)=\mathcal{F}\otimes_X^\L\mathcal{O}_X(n)
=\mathcal{F}
\otimes_X^\L
\overbrace{\mathcal{O}_X(1)\otimes^\L_X\cdots\otimes^\L_X\mathcal{O}_X(1)}^n.
\]

As we expect that the pair $(P_1,\rho)$ is ample,
reversing Serre's vanishing theorem,
we define  the following.

\begin{define}
The full subcategory $D^{\rho,\geq 0}$ 
\textup{(}resp.  $D^{\rho,\leq 0}$\textup{)} of
$D^b(\rmod)$ consists of objects $M^{\cdot}$ which satisfy
\begin{equation*}
\begin{split}
&\R\Hom^{\cdot}(P_1,\rho^nM^{\cdot})\in D^{\geq 0}(k\text{-}vect)
\quad
\text{for $n \gg 0$}\\
(\text{resp.}\quad
&\R\Hom^{\cdot}(P_1,\rho^nM^{\cdot})\in D^{\leq 0}(k\text{-}vect)
\quad
\text{for $n \gg 0$}).
\end{split}
\end{equation*}

\end{define}

\begin{rem}\label{above}
Since $A\cong P_1\oplus\rho P_1$,
$M^\cdot\in D^{\rho,\geq 0}\,
(\text{resp.}\,\,D^{\rho,\leq 0})$
if and only if
$\rho^nM^\cdot \in D^{\geq 0}(\rmod)\,
(\text{resp.}\,\,D^{\leq 0}(\rmod))$
for $n\gg 0$.
\end{rem}

\begin{prop}\label{t-st-prop}
The pair of full subcategories $(D^{\rho,\geq 0},D^{\rho,\leq 0})$
is a t-structure in $D^b(\rmod)$.
\end{prop}

\begin{proof}
Since the other conditions of t-structure
are obvious by Remark \ref{above},
it suffices to check the condition $(\text{t}3)$, namely,  
for every complex $M^{\cdot}\in D^b(\rmod)$
there is an exact triangle
\[
A^{\cdot}\longrightarrow M^{\cdot}\longrightarrow B^{\cdot}\shift
\]
such that $A^{\cdot}\in D^{\rho,\leq 0}$ and
$B^{\cdot}\in D^{\rho ,\geq 1}= D^{\rho, \geq 0}[-1]$.

By \cite[Lemma I.5.2]{Ha} 
an indecomposable object of $D^b(\rmod)$
is of the form $M[-i]$,
where $M\in \rmod$ and $i\in \Z$.
By Lemma \ref{pure},
$\rho^{n}M[-i]=M\otimes^{\L}_A\rho^{n}[-i]\in D^{[i-1,i]}(\rmod)$.
Hence
$M[-i]\in D^{\rho,\geq 1}$ for $i \geq 2$ and
$M[-i]\in D^{\rho,\leq 0}$ for $i \leq 0$.
The case when $i=-1$ is reduced to the following lemma.
\end{proof}

\begin{lemma}\label{t-st-lem}
For every $M\in \rmod$
there exists a submodule $M'$
such that
$\Hom(P_1,\rho^nM')=0$ for $n \gg 0$
and
$\Ext^{-1}_A(P_1,\rho^{n}M'')=0$ for $n\gg 0$,
where we set $M''= M/M'$.
\end{lemma}

\begin{proof}
Define
$\mathcal{T}_n=\{N\in \rmod\mid N\,\text{is generated by}\, \rho^{-n}S_0\}$
and $\mathcal{F}_n=\{N\in \rmod\mid \Hom_A(\rho^{-n}S_0,N)=0\}$.
Since
$\Ext^{1}_A(\rho^{-n}S_0,\rho^{-n}S_0)\cong\Ext^1_A(P_1,P_1)=0$
and $A$ is hereditary,
by \cite[Lemma III.4.2]{Ha},
the pair $(\mathcal{T}_n,\mathcal{F}_n)$ is a torsion theory
on the abelian category $\rmod$.
So,
if we set $t_n(M),f_n(M)$ to be
the image and cokernel of the canonical morphism
\[
\Hom(\rho^{-n}S_0,M)\otimes_k \rho^{-n}S_0\rightarrow M,
\]
then $t_n(M)\in \mathcal{T}_n,\,f_n(M)\in \mathcal{F}_n$.
By Lemma \ref{surj-lem} the canonical morphism 
\[
\Hom(\rho^{-(n+1)}S_0,\rho^{-n}S_0)\otimes_k\rho^{-(n+1)}S_0
\longrightarrow\rho^{-n}S_0
\]
is surjective,
hence $t_n(M)\subset t_{n+1}(M)$.
Since $\dim_k M < \infty$ ,
there is an integer  $n_0$ such that
$t_n(M)=t_{n_0}(M)$ for all $n\geq n_0$.
Then by definition $f_n(M)=f_{n_0}(M)$ for all $n \geq n_0$.
Set $M'=t_{n_0}(M),M''=f_{n_0}(M)$.
Then $M'\in\mathcal{T}_{n},M''\in\mathcal{F}_{n}$ for all $n\geq n_0$.

It is easy to see that
$\Hom(P_1,\rho^n N)\cong\Ext^1(\rho^{-(n-1)}S_0,N)=0$ 
for $N\in \mathcal{T}_{n-1}$
and 
$\Ext^{-1}_A(P_1,\rho^nN)\cong\Hom_A(\rho^{-(n-1)}S_0,N)=0$
for $N\in \mathcal{F}_{n-1}$.
Therefore 
$\Hom(P_0,\rho^nM')=0,\Ext^{-1}_A(P_0,\rho^{n}M'')=0$ 
for all $n > n_0$.
\end{proof}

\begin{rem}
Set
$\mathcal{T'}=\{N\in\rmod\mid\Hom(P_1,\rho^nN)=0\,\text{for}\,\,n \gg 0\}$
and
$\mathcal{F'}=\{N\in\rmod\mid\Ext^{-1}(P_1,\rho^nN)=0\,\text{for}\,\,n \gg 
0\}$.
Then by the Lemma \ref{t-st-lem} it is easy to see that
$(\mathcal{T'},\mathcal{F'})$ is a torsion pair on $\rmod$.

From this torsion pair 
we can define a t-structure in $D^b(\rmod)$ 
by setting 
\begin{equation*}
\begin{split}
D'^{\geq 0}:=&
\{M^\cdot\in D^{\geq 0}(\rmod)\mid \H^0(M^\cdot)\in\mathcal{F'}\}\\
D'^{\leq 0}:=&
\{M^\cdot\in D^{\leq 1}(\rmod)\mid \H^1(M^\cdot)\in\mathcal{T'}\}.
\end{split}
\end{equation*}
\textup{(}See \textup{\cite[Proposition I.2.1]{HRS}}\textup{)}.
However, this is not a new t-structure.
It is easy to see that
$(D'^{\geq 0},D'^{\leq 0})=(D^{\rho,\geq 0},D^{\rho, \leq 0})$.
\end{rem}

Let $\mathcal{H}^\rho$ be the heart of the t-structure  
$(D^{\rho,\geq 0},D^{\rho,\leq 0})$.
By Proposition \ref{pure}, $P_0,P_1\in \mathcal{H}^{\rho}$.
Now let us  consider 
the triple $(\mathcal{H}^\rho,P_1,\rho)$.

\begin{prop}\label{rho-ample}
The pair $(P_1,\rho)$ is ample on 
$(\mathcal{H}^\rho,P_1,\rho)$.
\end{prop}

\begin{proof}
We check the conditions (1) and (2) of Definition \ref{ample}.

First note that the cokernel of the morphism
$f: M^\cdot\longrightarrow N^\cdot$ in the abelian category
$\mathcal{H}^\rho$ is
$\tau^{\rho}_{\geq 0}(\Cone(f))$,
where
$\tau^{\rho}_{\geq 0}: D^b(\rmod)\longrightarrow D^{\rho,\geq 0}$
 is the truncation functor
(See \cite[IV.4]{GM}).
So $f$ is surjective in $\mathcal{H}^{\rho}$
if and only if
$\Cone(f)\in D^{\rho,\leq -1}$.

$(1)$
It is easy to see that for every $M\in\rmod$
there is an exact sequence
\[
0\longrightarrow\Hom_A(P_1,M)\otimes_k P_1
\longrightarrow M
\longrightarrow\Hom_A(P_0,M)\otimes_k S_0
\longrightarrow 0
\]
which is functorial in M.
So for every $M^{\cdot}\in D^b(A\text{-mod})$ and for every $n$
there is an exact triangle
\[
\R\Hom(P_1,\rho^nM^\cdot)\otimes_k\rho^{-n}P_1\longrightarrow
M^\cdot \longrightarrow
\R\Hom(P_0,\rho^{n-1}M^\cdot)\otimes_k\rho^{-n}S_0\shift .
\]

Since
$\rho^{-n}P_1\in D^{\rho,\leq n}$ for every integer $n$,
$\rho^{-n}S_0\cong\rho^{-(n+1)}P_1[1]\in D^{\rho,\leq -1}$.

Now assume that $M^\cdot\in\mathcal{H}^{\rho}$
and take an integer $n$ such that the complex
$\R\Hom(P_1,\rho^{n-1}M^{\cdot})$ is pure.
Then
$\R\Hom(P_1,\rho^{n-1}M^\cdot)\otimes_k\rho^{-n}S_0\in D^{\rho,\leq -1}$.

$(2)$
Let
\[
M^{\cdot}\stackrel{f}{\longrightarrow}
N^{\cdot}\longrightarrow
L^{\cdot}\shift
\]
be a triangle
such that
$M^{\cdot},N^{\cdot}\in\mathcal{H}^{\rho}$ and
$L^{\cdot}\in D^{\rho,\leq -1}$.
Take an integer $n$ 
such that $\Hom(P_1,\rho^nL^{\cdot})=0$.
Then the induced morphism
$\Hom(P_1,\rho^nM^{\cdot})\longrightarrow
\Hom(P_1,\rho^nN^{\cdot})$
is surjective.
\end{proof}

Combining Proposition \ref{coordinate} and Proposition \ref{rho-ample}
and applying Theorem \ref{Serre} we obtain the following proposition.

\begin{prop}
The triple 
$(\mathcal{H}^\rho,P_1,\rho)$ is equivalent
to the noncommutative projective scheme 
$\proj R=(\qgr R,\overline{R},(1))$.
i.e.,
the functor 
\[
\overline{\Gamma}_*:\mathcal{H}^{\rho}
\longrightarrow\qgr R,\quad
\overline{\Gamma}_*(M^\cdot)=\oplus_{n\geq 0}\R\Hom(P_1,\rho^nM^{\cdot})
\]
is an equivalence of categories,
$\overline{\Gamma}_*(P_1)=\overline{R}$ and 
$\overline{\Gamma}_*\circ\rho=(1)\circ\overline{\Gamma}_*$.
\end{prop}

Now we complete the proof of our main theorem.

The functor
\[
\Gamma_* : \rmod \longrightarrow \Gr R ,\quad
\Gamma_*(M):= \bigoplus_{n\geq 0}\Hom(P_1,M\otimes_A \rho^n)
\]
is right exact. Let $\L\Gamma_*$ be the left derived functor of $\Gamma_*$.
A finite projective $A$-module $P$ is of the form
$P\cong P_1^{\oplus n}\oplus P_0^{\oplus m}$ for some $n,m\in\Z_{\geq 0}$.
Then $\Gamma_*(P)\cong R^{\oplus n}\oplus R(1)^{\oplus m}$ is a
finite $R$-module.
Since every complex $M^\cdot\in D^b(\rmod)$ is isomorphic to 
a bounded complex of finitely generated projective $A$-modules,
$\L\Gamma_*(M^\cdot)$ 
has cohomologies which are finitely presented.
Thus there is an exact functor
\[
\L\Gamma_* : D^b(\rmod)\longrightarrow D^b_{\gr R}(\Gr R),\quad
\L\Gamma_{*}(M^\cdot)=\bigoplus_{n\geq 0}\Hom(P_1,\rho^{n}M^\cdot).
\]

Since 
the quotient functor $\pi: \Gr R\longrightarrow \QGr R$
is exact by \cite[Theorem 4.3.8]{Pop},
this functor can be extended  to
the exact functor
$\pi: D^b_{\gr R}(\Gr R)\longrightarrow D^b_{\qgr R}(\QGr R)$.
Define $\L\overline{\Gamma}_*:=\pi\circ\L\Gamma_*$.

The canonical functor
$c : D^b(\cohproj R)\longrightarrow D^b_{\cohproj R}(\QGr R)$
is an equivalence by \cite[Lemma 4.3.3]{BvB}.
Hence, to prove our main theorem,
we have only to show that 
$\L\overline{\Gamma}_*$ is an equivalence.

There is a commutative diagram 
\begin{equation}\label{sikaku}
\begin{xymatrix}
{
D^b(\rmod)\ar[rr]^{\L\overline{\Gamma}_*}
&&D^b_{\cohproj R}(\QGr R)\\
\mathcal{H}^\rho \ar@{^{(}->}[u]^{i_A} \ar[r]^{\sim}_{\overline{\Gamma}_*} &
\cohproj R\, \ar@{^{(}->}[r]_{i_R} & 
D^b(\cohproj R)\ar[u]^c_{\wr}
}
\end{xymatrix}
\end{equation}
where $i_A,i_R$ are inclusions.

It is clear that $\cohproj R$ generates $D^b(\cohproj R)$,
i.e.,
the minimal triangulated full subcategory of $D^b(\qgr R)$ 
containing $\qgr R$ is $D^b(\qgr R)$ itself.
Therefore 
the image 
$(\L\overline{\Gamma}_*\circ i_A)(\mathcal{H}^\rho)$ 
generates 
$D^b_{\qgr R}(\QGr R)$.
Thus $\L\overline{\Gamma}_*$ is essentially surjective.
To prove $\L\overline{\Gamma}_*$ is an equivalence 
it suffices to show that 
$\L\overline{\Gamma}_*$ is fully faithful.

Since every complex $M^\cdot\in D^b(\rmod)$ is 
obtained from $P_0,P_1$ 
by taking finite number of cones and shifts, 
the problem is reduced to the following lemma.

\begin{lemma}\label{isom-lem}
The map 
\[
\Hom_{D^b(\rmod)}(P_i,P_j[l])
\xrightarrow{\L{{\overline{\Gamma}}_{*}}_{P_i,P_j[l]}} 
\Hom_{D^b_{\qgr R}(\QGr R)}
(\L\overline{\Gamma}_{*}(P_i),\L\overline{\Gamma}_{*}(P_j)[l])
\]
is an isomorphism for every $i,j\in\{0,1\}$ and $l\in \Z$.
\end{lemma}

\begin{proof}
For the case $l=0$,
by the commutativity of the diagram (\ref{sikaku}),
the map $\L{\overline{\Gamma}_{*}}_{P_i,P_j}$ is 
equal to 
\begin{equation*}  
\begin{split}
\Hom_{D^b(\rmod)}(P_i,P_j)&\cong 
\Hom_{\mathcal{H^\rho}}(P_i,P_j)\\
&\cong
\Hom_{\qgr R}(\overline{\Gamma}_{*}(P_i),\overline{\Gamma}_{*}(P_j))\\
&\cong
\Hom_{D^b_{\qgr R}(\QGr R)}
(\L\overline{\Gamma}_{*}(P_i),\L\overline{\Gamma}_{*}(P_j)).
\end{split}
\end{equation*}
Hence $\L{{\overline{\Gamma}}_{*}}_{P_i,P_j}$ is an isomorphism.
  
By \cite[Appendix]{Be2}
there exists an exact functor
$F : D^b(\mathcal{H}^\rho)\longrightarrow D^b(\rmod)$
and a commutative diagram
\[
\begin{xymatrix}
{
D^b(\mathcal{H}^\rho)\ar[r]^F&D^b(\rmod)\\
\mathcal{H}^\rho\ar@{^{(}->}[u]^{i_H} \ar@{^{(}->}[ur]^{i_A}
}
\end{xymatrix},
\]
where $i_H,i_A$ are inclusions.

The restriction $\L\overline{\Gamma}_{*}\circ F|_{\mathcal{H}^\rho}$
to $\mathcal{H}^\rho$ is equal to
$
c\circ i_R \circ \overline{\Gamma}_{*}
$ 
with the notation in the diagram (\ref{sikaku}).
Hence 
$\L\overline{\Gamma}_{*}\circ F$ is an equivalence.
Therefore, the map 
\begin{equation*}
\Hom_{D^b(\rmod)}(M^\cdot,N^\cdot[l])
\xrightarrow{\L{{\overline{\Gamma}_{*}}}_{M^\cdot,N^\cdot[l]}}  
\Hom_{D^b_{\qgr R}(\QGr R)}
(\L\overline{\Gamma}_*(M^\cdot),\L\overline{\Gamma}_*(N^\cdot)[l])
\end{equation*}
is surjective for $M^\cdot,N^\cdot\in D^b(\rmod)$
and $l\in\Z$.
Now it is easy to prove the lemma 
because 
$\Hom_{\rmod}(P_i,P_j[l])=0$ 
for every $i,j\in\{0,1\}$ and $l\neq 0$.
\end{proof}

\begin{theorem}[Theorem \ref{main}]\label{main2}
Let $k$ be an algebraically closed field,
and  $A=k\overrightarrow{\Omega}_N$,  the path algebra of
the $N$-Kronecker quiver $\overrightarrow{\Omega}_N$.
 Set  moreover $R=k\langle X_1,\dots,X_N \rangle/(\sum^N_{i=1}X_i^2)$.
 Then, for $N\geq 2$, the functor
\[
\L\overline{\Gamma}_{*}:D^b(\rmod)\longrightarrow D^b_{\qgr R}(\QGr R)
\]
is an equivalence of triangulated categories.
\end{theorem}

\begin{rem}
From the argument above,  
$\{\overline{R},\overline{R(1)}\}$ 
is a full strongly exceptional sequence,
and by  \cite[Theorem 4.5.2]{Pop}
the abelian category $\QGr R$ 
has enough injectives.
Hence we can apply  
\textup{\cite[Theorem 6.2]{Bo}}
\textup{(}and the following Remark
\textup{)}
and prove that
the functor 
\[
\R\Hom(T,-):D^b(\qgr R)\longrightarrow D^b(\rmod)
\]
is an equivalence of triangulated categories 
where we set $T=\overline{R}\oplus\overline{R(1)}$.

Since $A\cong \End(T)$,
$T$ has a natural left $A$-module structure. 
It is known that 
the functor
$-\otimes^\L_AT$ is a quasi-inverse to the functor
$\R\Hom(T,- )$.

It is easy to see that there is a  commutative diagram
\begin{equation*}
\begin{xymatrix}
{
D^b(\rmod)\ar[r]^{\L\overline{\Gamma}_{*}} \ar[dr]_{-\otimes^\L_A T}
&D^b_{\cohproj R}(\QGr R)\\
 &D^b(\cohproj R).\ar[u]^c_{\wr}
}
\end{xymatrix}
\end{equation*}
Thus the functors $\L\overline{\Gamma}_{*}$ and $-\otimes^\L_AT$
are essentially isomorphic.
\end{rem}

{\scshape
\begin{flushright}
\begin{tabular}{l}
Department of Mathematics, \\
Kyoto University, Kyoto 602-8502, \\
Japan \\
{\upshape e-mail: minamoto@kusm.kyoto-u.ac.jp}
\end{tabular}
\end{flushright}
}


\begin{thebibliography}{Sch1}

\bibitem[AZ]{NPS}
M. Artin, and J.J. Zhang,
Noncommutative projective schemes,
Advances in Math.
\textbf{109} (1994),
pp. 228-287.



\bibitem[Be]{Be}
A.A. Beilinson,
Coherent sheaves on $\P^n$ and problems of linear algebras,
Func. Anal. Appl. \textbf{12} (1978),
pp. 214-216.

\bibitem[Be2]{Be2}
A.A. Beilinson,
On the derived category of perverse sheaves,
in K-theory, arithmetic and geometry (Moscow 1984-1986),
pp. 27-41,
Lecture Notes in Math. \textbf{1289},
Springer, Berlin, 1987.

\bibitem[Bo]{Bo}
A.I. Bondal,
Representations of associative algebras
and coherent sheaves,
Izv. Akad. Nauk SSSR Ser. Mat. \textbf{53} (1989),
25-44.


\bibitem[BvB]{BvB}
A. Bondal, and M. Van den Bergh,
Generators and representability of functors in commutative
and noncommutative geometry,
Mosc. Math. J.,
\textbf{3} (2003),
no. 1,
pp. 1-36.

\bibitem[GM]{GM}
S.I. Gelfand, and Yu.I. Manin,
Methods of Homological Algebras,
(Second edition)
Springer Monographs in Mathematica.
Springer-Verlag, Berlin,2003.


\bibitem[Hap]{Ha}
D. Happel,
Triangulated Categories in the Representation Theory of
Finite-Dimensional Algebras,
London Math. Soc. Lecture Notes \textbf{119},
University Press,Cambridge, 1987.

\bibitem[HRS]{HRS}
D. Happel, I. Reiten, and S. Smalo,
Tilting in abelian categories and quasitilted algebras,
Memoirs of the AMS, vol.\textbf{575},
Amer. Math. Soc., 1996.

\bibitem[Har]{AG}
R. Hartshorne,
Algebraic Geometry,
Graduated Texts in Mathematics,
Springer-Verlag,1977.

\bibitem[KR]{KR}
M. Kontsevich, and A. Rosenberg,
Noncommutative smooth spaces,
The Gelfand Seminars 1996-1999,pp. 85-108,
Gelfand Math Sem.,
Birkha$\ddot{u}$ser Boston,
Boston, MA, 2000.

\bibitem[MY]{MY}
J. Miyachi, and A. Yekutieli,
Derived Picard groups of finite dimensional hereditary algebras,
Compositio Math. \textbf{129} (2001),
no. 3, pp. 341-368.

\bibitem[O]{O}
D. Orlov,
Derived categories of coherent sheaves and triangulated categories of singularities,
math.AG/0503632.


\bibitem[Po]{Po}
A. Polishchuk,
Nonncommutative proj and coherent algebras,
Math. Res. Lett.,
\textbf{12} (2005),1,
pp. 63-74.

\bibitem[Pop]{Pop}
N. Popescu,
Abelian Categories with Applications to Rings and Modules,
L.M.S. Monographs, Academic Press,
New York, 1973.

\bibitem[SvB]{SvB}
J.T. Stafford, and M. Van den Bergh,
Noncommutative curves and noncommutative surfaces,
Bull. AMS \textbf{38} (2001), pp. 171-216.

\bibitem[Y]{Y}
A. Yekutieli,
Dualizing complexes, Morita equivalence and the derived Picard group
of a ring,
J. London Math Soc.
\textbf{60} (1999), pp. 723-746.
  
\bibitem[Z]{Z}
J. J. Zhang,
Twisted graded algebras and equivalences of graded categories,
Proc. London math. soc.
(3)\textbf{72}(1996), no,2.
pp. 281-311

\end{thebibliography}
\end{document}